\documentclass[12pt]{amsart}
\usepackage{amsmath,amsfonts,amssymb,amscd,amsthm,amsbsy,epsf}
\newtheorem*{thm}{Theorem} % [subsection]
\newtheorem*{cor}{Corollary} % [subsection]
\newtheorem*{Lem}{Lemma}
\def\Int{\operatorname{int}}
\def\R{{\mathbb  R}}
\def\zed{{\mathbb Z}}
\begin{document}

\centerline{\bf ON EMBEDDING INFINITE CYCLIC COVERS}
\smallskip
\centerline{\bf IN COMPACT 3-MANIFOLDS}
%On embedding infinite cyclic covers in compact 3-manifolds}
\bigskip
\centerline{\sc C. McA. Gordon}

\bigskip

\subsection*{1.}		% 1.

By a knot we shall mean a smooth knot $K$ in $S^3$. 
The exterior of $K$ is $E(K) = S^3 - \Int N(K)$, where $N(K)$ is a tubular 
neighborhood of $K$, and we shall refer to the infinite cyclic cover of $E(K)$ 
as simply the infinite cyclic cover of $K$.

In [JNWZ] the authors consider the question (attributed to John Stallings) of 
when the infinite cyclic cover of a knot embeds in $S^3$. 
Following [JNWZ], we say that such a knot has Property IE. 
Clearly fibered knots have Property IE. 
In [JNWZ] it is shown that if a non-fibered genus 1 knot has Property IE then 
its Alexander polynomial is either 1 or $2-5t+2t^2$, and that in each case 
there are infinitely many genus~1 knots having Property IE.

[JNWZ] also raises the question of whether the infinite cyclic cover of a 
knot always embeds in some compact 3-manifold. 
The following theorem provides a negative answer.

\begin{thm} %Theorem. 
The infinite cyclic cover of the untwisted Whitehead double of a non-trivial 
knot does not embed in any compact $3$-manifold.
\end{thm}

\begin{cor} %Corollary. 
There are infinitely many genus $1$ knots with Alexander polynomial $1$ 
whose infinite cyclic covers do not embed in any compact $3$-manifold.
\end{cor}

\subsection*{2}	%2.

Let $L = J \cup J'$ be the Whitehead link, and $W$ be its exterior; see 
Figure~1. 
Let $T$ be the boundary component of $W$ corresponding to the component $J$ 
of $L$, and let $m,\ell$ be a meridian-longitude pair for $J$ on $T$. 
Let $k$ be a non-trivial knot, with exterior $X$, and let $\mu,\lambda$
be a meridian-longitude pair for $k$ on  $\partial X$. 
Let $K$ be the untwisted Whitehead double of $k$. 
Then the exterior of $K$ is $Y = W \cup_T  X$, where $T$ is identified with  
$\partial X$ so that $m\leftrightarrow \lambda$, $\ell\leftrightarrow \mu$.
$$\begin{matrix}\epsfysize=1.7truein\epsfbox{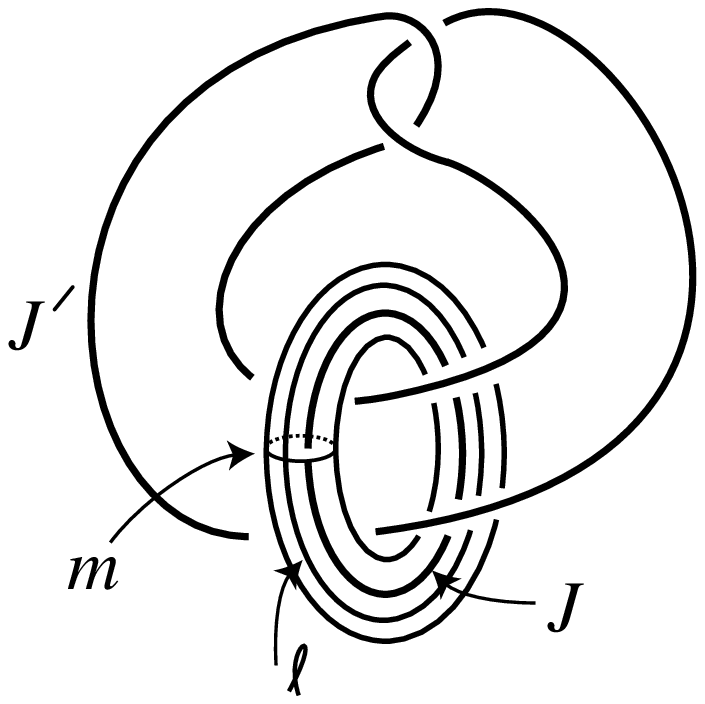}\\
\text{Figure 1}\end{matrix}$$

The exterior of $J$ (resp. $J'$) in $S^3$ is a solid torus $V$ (resp. $V'$). 
Let $V'_\infty$ be the infinite cyclic cover of $V'$, and let $\tilde J$ 
be the inverse image of $J$ in $V'_\infty$. 
The infinite cyclic cover of $J'$ in the solid torus $V$ is 
$W_\infty  = V'_\infty - \Int N(\tilde J)$. 
Since there is an isotopy of $S^3$ interchanging the components of $L$, the 
link $\tilde J$ in $V'_\infty \cong  D^2 \times\R$ is as shown in 
Figure~2. 
$$\begin{matrix}\epsfysize=3.5truein\epsfbox{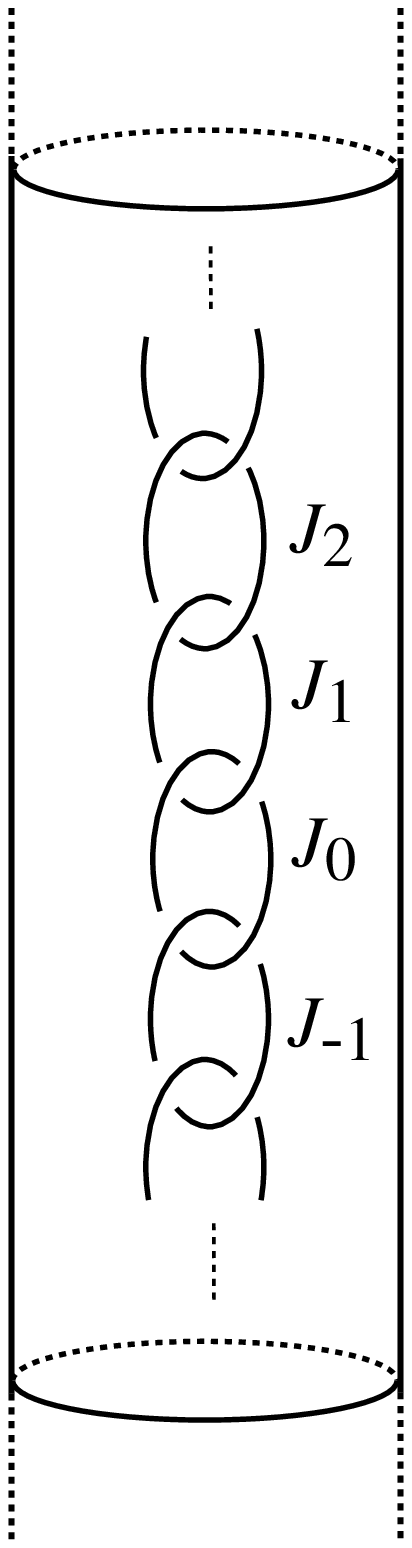}\\
\text{Figure 2}\end{matrix}\hskip1truein 
\begin{matrix} \\ \noalign{\vskip18pt}
\epsfysize=2.5truein\epsfbox{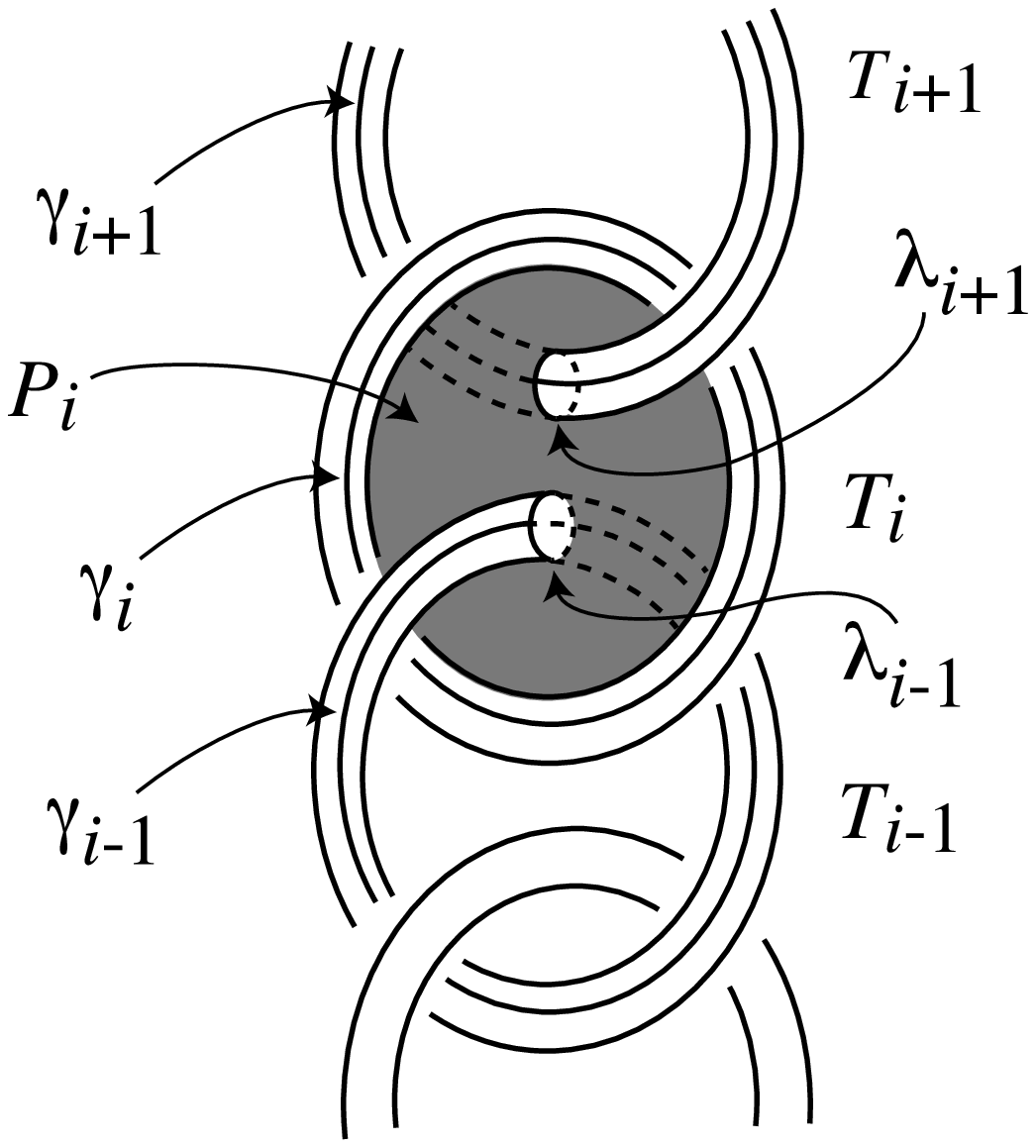}\\
\noalign{\vskip38pt}  
\text{Figure 3}\end{matrix}$$

Let the components of $\tilde J$ be $J_i$, $i \in\zed$, 
indexed in the obvious way, and let $T_i =\partial N(J_i)$. 
Then the infinite cyclic cover of $Y$ is 
$Y_\infty  = W_\infty \cup (\coprod_{i\in\zed} X_i)$, 
where $X_i$ is a copy of $X$, attached to $W_\infty$  
by identifying $T_i$ with   $\partial X_i$ via
$\tilde m_i\leftrightarrow \lambda_i$, $\tilde \ell_i\leftrightarrow \mu_i$,
where $\mu_i,\lambda_i$ is a meridian-longitude pair on $\partial  X_i$, and
$\tilde m_i,\tilde\ell_i$ are lifts of $m,\ell$. 
Note that $\tilde m_i$   is a meridian of $J_i$.

Let $\gamma_i\subset  T_i$ be a longitude of $J_i$ in $D^2\times\R$. 
Then there is a planar surface $P_i\subset W_\infty$ with $\partial  P_i =
\gamma_i\cup\lambda_{i-1} \cup\lambda_{i+1}$.
See Figure~3. 

Let $S_j \subset X_j$ be a copy of a Seifert surface for $k$, and 
define $F_i = P_i \cup S_{i-1}\cup S_{i+1}$. 
Then $F_i$ is an orientable surface in $Y_\infty$ with 
$F_i\cap X_i = \partial  F_i = \gamma_i$. 
Note that each of $\gamma_{i-1}$ and $\gamma_{i+1}$ intersects $F_i$ 
transversely in a single point.

We are now ready to prove the theorem. 
So suppose that $Y_\infty$ embeds in a compact 3-manifold $M$. 
By passing to a double cover if necessary we may assume that $M$
is orientable. 

\begin{Lem}  % Lemma 2. 
$T_i$ is incompressible in $M$.
\end{Lem}

\begin{proof} 
Let $D$ be a compressing disk for $T_i$ in $M$. Since $k$ is non-trivial,
$D\cap  X_i = \partial  D$. 
Recall the orientable surface $F_i \subset Y_\infty$ 
with $\partial F_i \cap  X_i =
\partial F_i = \gamma_i$. 
Then $\partial D$ and $\gamma_i$ are both null-homologous in 
$\overline{M - X_i}$, 
and hence they are isotopic on $T_i$. 
Let $\widehat F_i$ be the (singular) closed surface $F_i\cup D$. 
We may assume that $D$ intersects $T_{i+1}$ transversely in a finite number 
of simple closed curves, each essential on $T_{i+1}$.

If $D\cap  T_{i+1} = \emptyset$, then $\gamma_{i+1}$ intersects $\widehat F_i$ 
transversely in a single point, and so $[\gamma_{i+1}]$ has infinite order 
in $H_1(M)$. 
%But if $T_{i+1}$ were compressible, then, as argued above, $\gamma_{i+1}$ 
%would bound a disk in $M$, a contradiction.
But this contradicts the fact that $\gamma_{i+1}$ bounds $F_{i+1}$.

If $D\cap  T_{i+1} \ne\emptyset$, then an innermost disk on $D$ is a 
compressing disk $D'$ for $T_{i+1}$ disjoint from $T_i$. 
Now, $\gamma_i$ meets $\widehat F_{i+1} = F_{i+1} \cup D'$
transversely in a single point, and we get a contradiction as before.
\end{proof}

\begin{proof}[Proof of Theorem] 
By the Lemma and Haken Finiteness 
there exist distinct $i,j,k$ such that $T_i$, $T_j$ and $T_k$ are 
mutually parallel in $M$. 
By relabeling if necessary we may assume that $T_j$ is contained in the 
interior of the product region between $T_i$ and $T_k$. 
Since $T_j$ bounds $X_j$, either $T_i$ or $T_k$ is contained in $X_j$, 
which is absurd.
\end{proof}

\bigskip
\bigskip

{\footnotesize\parindent=0pt
Department of Mathematics \par
The University of Texas at Austin \par
1 University Station C1200 \par
Austin TX 78712-0257  USA\par
\medskip

gordon@math.utexas.edu}

\end{document}